\documentclass[final]{siamart171218}
\usepackage{moreverb,amsmath, amsfonts, bm, bbm, graphicx, caption, subcaption,pgf}

\begin{document}

\title{Scalable Extended Dynamic Mode Decomposition using Random Kernel Approximation}
\author{Anthony M. DeGennaro%
        \thanks{Computational Science Initiative, Brookhaven National Laboratory (\email{adegennaro@bnl.gov}).}%
        \and
        Nathan M. Urban%
        \thanks{Computer, Computational, and Statistical Sciences, Los Alamos National Laboratory (\email{nurban@lanl.gov}).}}

\date{\today}

\newcommand*{\vertbar}{\rule[-1ex]{0.5pt}{2.5ex}}
\newcommand*{\horzbar}{\rule[.5ex]{2.5ex}{0.5pt}}
\newcommand{\amd}[1]{\textcolor{green}{#1}}
\newcommand{\old}[1]{\textcolor{red}{#1}}

\maketitle

\begin{abstract}
The Koopman operator is a linear, infinite-dimensional operator that governs the dynamics of system observables; Extended Dynamic Mode Decomposition (EDMD) is a data-driven method for approximating the Koopman operator using functions (features) of the system state snapshots. This paper investigates an approach to EDMD in which the features used provide random
approximations to a particular kernel function. The objective of this
is computational economy for large data sets: EDMD is generally
ill-suited for problems with large state dimension, and its dual kernel formulation (KDMD) is well-suited for
such problems only if the number of data snapshots is relatively
small. We discuss two specific methods for generating features: random Fourier
features, and the N{\"y}strom method. The first method is a
data-independent method for translation-invariant kernels only and
involves random sampling in feature space; the second method is a
data-dependent empirical method that may be used for any kernel and
involves random sampling of data. We first discuss how these ideas may be applied in an EDMD context, as well as a means for adaptively adding
random Fourier features. We demonstrate these methods on two
example problems and conclude with an analysis of the relative benefits and drawbacks of each method.
\end{abstract}

\section{Introduction}
Dynamic Mode Decomposition (DMD)~\cite{Schmidt2008,SchmidJFM,Tu2013_journal,chen} is a
popular technique in model reduction for nonlinear dynamical
systems. The goal of this method is to produce an approximation of the
Koopman operator~\cite{Rowley2009,mezic2013,mezic2005,koopman} -- a linear, infinite-dimensional
operator that propagates system observables -- using temporal
snapshots of the system observables. As a linear operator, the Koopman
operator may be completely described in terms of its eigenspectrum,
and so the DMD process can be thought of as a data-driven method to
approximate the leading eigenmodes, eigenvalues, and
eigenfunctions. If the leading eigenfunctions of Koopman lie within
the span of the observable state data used by DMD, then DMD can
produce a good approximation of them. However, as pointed out
in~\cite{WilliamsEDMD_journal}, DMD is only capable of producing an
approximation of the eigenfunctions using linear monomials in the
state. While this may be sufficient in cases where linear behavior
dominates the dynamics, in many cases it is not, and so more
sophisticated techniques must be used to approximate Koopman.

One such
technique approaches this issue by generalizing the DMD machinery so
that one is not limited to using system state data, but in fact may
use functions (i.e., features) of the state data that may be more
well-suited to approximating the Koopman eigenfunctions. This method
is known as Extended DMD (EDMD)~\cite{WilliamsEDMD_journal}. While certainly
helpful in many contexts, EDMD is fundamentally a spectral method and
hence suffers from the ``curse of dimensionality'': the computational
expense is dominated by the number of basis functions
(features), which tends to grow quickly with increasing state
dimension. This problem was partially solved by utilizing the
so-called ``kernel-trick'', which led to Kernel DMD
(KDMD)\cite{WilliamsKDMD_journal}, in which one need not explicitly transform
the data into feature space; all that is required is knowledge of a
kernel function that operates on the state data to produce inner
products in feature space. However, while the computational expense of
KDMD is not dominated by the number of basis functions, it is
dominated by the number of snapshots and the state size. Hence,
neither EDMD or KDMD are particularly well-suited to scenarios in
which the state dimension and number of snapshots are both large.

The objective of this work is to make progress toward mitigating the
computational expense of EDMD for scenarios in which the state size
and number of snapshots are large. We propose doing this by using
random features that generate efficient approximations to a particular
kernel function. We will discuss this in greater detail later, but
briefly, this choice is motivated by the theorems of Mercer and
Bochner, which guarantee spectral decompositions of kernel functions
under certain conditions. In choosing to use a random collection of
features that approximate a particular kernel as the features in an
EDMD method, we are conceptually producing a Monte Carlo approximation
of the KDMD method, using EDMD. The intended advantage of this is,
loosely, a blending of the best features of EDMD and KDMD: we wish to
use the EDMD method so that we are not algorithmically limited by the
number of snapshots, but also desire an efficent EDMD basis that does
not scale as badly with state dimension as some of the more
traditional feature choices (e.g., polynomial tensor products or
radial basis functions (RBFs)). This claim is based on the hope that
the chosen kernel function may be efficiently approximated with a
modest number of features, and that therefore, the number of
EDMD features required for a faithful Koopman approximation is modest
as well. It is reasonable to expect that both the state size and
complexity of the system dynamics could affect the rate of convergence
of the EDMD algorithm; however, as we will see in the numerical
examples, EDMD using kernel approximation methods can produce
efficient, accurate Koopman approximations in benchmark problems with
a state dimension large enough that most feature choices -- like
polynomials or RBFs -- would be computationally intractable.

Depending on our choice of kernel, there are potentially two
approaches for generating random features for EDMD: random Fourier
features~\cite{Rahimi2007,Rahimi2008}, and the N{\"y}strom
method~\cite{Nystrom1,Nystrom2}. In both methods, one begins by
choosing a desired kernel function (which sets inner products in
feature space). If the chosen kernel is positive semi-definite,
symmetric, and translation-invariant, then Bochner's theorem and the
method of random Fourier features prescribes a representation in the
Fourier basis, with the modes adhering to a distribution related to
the kernel. In such a case, the EDMD features would consist of a
finite collection of those modes, drawn randomly from their underlying
distribution. If, instead, our kernel is not translation-invariant
(but is positive semi-definite and symmetric), then Mercer's theorem
guarantees the existence of an eigendecomposition, which we may
estimate empirically from data using the N{\"y}strom method. In this
more general case, the EDMD features would consist of the approximate
kernel eigenfunctions.

To state the basic problem more concretely for background purposes, we assume we have access to a time-series of $(M+1)$ state snapshot vectors $\lbrace x_0 , x_1, \dots , x_M \rbrace$, where any snapshot $x_i \in \mathbb{R}^d$. These data may have been generated by a nonlinear dynamical process, but DMD (and its variants) assumes that there is a corresponding linear dynamical system that can approximate those dynamics. In the case of DMD, this system is assumed linear in the state:
\begin{equation}
Y = AX
\end{equation}
where $X = \left[ x_0 , x_1 , \dots , x_{M-1} \right] \in \mathbb{R}^{d \times M}$ and $Y = \left[ x_1 , x_2 , \dots , x_M \right] \in \mathbb{R}^{d \times M}$ are the snapshot matrices and $A \in \mathbb{R}^{d \times d}$ is the DMD matrix. In the case of EDMD, following Williams~\cite{WilliamsEDMD_journal}, the system is assumed to be linear in some chosen feature space. That is, given a dictionary of $K$ functions acting on the state $\lbrace \psi_1(x) \dots \psi_K(x) \rbrace$, where $\psi_i(x) : \mathbb{R}^d \mapsto \mathbb{R}$, we define the vector valued function $\Psi : \mathbb{R}^{d} \mapsto \mathbb{R}^{1 \times K}$ as $\Psi(x) = \left[ \psi_1(x) \dots \psi_K(x) \right]$ and assume the system dynamics are linear when written as:
\begin{equation}
\label{eq:EDMDproblem}
\Psi_Y = \Psi_XA
\end{equation}
where $\Psi_X,\Psi_Y\in \mathbb{R}^{M \times K}$ are the feature matrices, with row $j$ corresponding to $\Psi(x_j)$, and $A \in \mathbb{R}^{K \times K}$ is the EDMD matrix. The EDMD Koopman matrix can thus be calculated as:
\begin{equation}
\label{eq:EDMDKoopman}
A = \Psi_X^\dagger \Psi_Y
\end{equation}
The computation of $A$ in Eq.~\ref{eq:EDMDKoopman} is easier by means
of the equivalent formulation:
\begin{equation}
\label{eq:GramKoopman}
A = G^\dagger H
\end{equation}
where $G = \Psi_X^T \Psi_X$ and $H = \Psi_X^T \Psi_Y$. As $G,H \in \mathbb{R}^{K \times K}$, the cost of
Eq.~\ref{eq:GramKoopman} is determined by $K$, and so the approaches we
discuss in this paper -- random Fourier features and the N{\"y}strom
approximation -- are aimed at producing a good approximation to $A$
using an economical number $K$ of features.

To be clear -- we are not attempting to improve the asymptotic scaling
trends of EDMD with $K$, $M$, or $d$. Our goal is simply to
apply a set of methods from machine learning in order to generate a
basis which requires a relatively low value of $K$ for problems where
both $d$ and $M$ are large.

We briefly pause here in order to step back and frame our work within
a larger context. One of the reasons that Koopman analysis shows such
promise in dynamical systems modeling is that it is a framework that
can be approached in a purely data-driven manner, and data-driven
learning is a paradigm that has enjoyed great progress
recently. Indeed, advances in machine learning have demonstrated the
power of empirically-trained models in a wide-range of computer
science and applied mathematics applications, including computer
vision~\cite{krizhevsky,simonyan,girshick}, speech
recognition~\cite{hinton,dahl}, manifold
learning~\cite{coifman,lee,nadler,nadler2}, and stochastic process
modeling~\cite{rasmussen,mackay}. It is natural to hope to import the
machinery underpinning these successes into the Koopman setting, and
our work is a part of that effort.

Recent research has produced many computational approaches to
approximating Koopman, including Generalized Laplace
Analysis~\cite{budisic,mauroy,mauroy2}, the Ulam Galerkin
Method~\cite{froyland,bollt}, and DMD. Of these, DMD is arguably the
most popular, owing to its ease of implementation and the fact that it
is well-suited to machine learning extensions (e.g.,
see~\cite{WilliamsEDMD_journal,WilliamsKDMD_journal,otto}). While
these properties make DMD a promising methodology, overwhelmingly
large amounts of data can be problematic. For this reason, there is a
community of researchers interested in developing scalable algorithms
for data-driven tasks. In point of fact, a number of recent studies
attempt to alleviate the computational burden of DMD, by leveraging
the general concepts of randomization, compressive sampling, and/or
sparsity; in particular,
see~\cite{ErichsonRandomizedDMD2017,ErichsonRandomizedDMD2016,ErichsonCompressedDMD,BruntonCompressiveDMD,gueniat2015,jovanovic,wynn,Hemati2014}. Of
course, there are other aspects of DMD that require improvement
besides scalability (see, for example, recent work in making DMD
robust to noise~\cite{dawson,hemati2017,bagheri}). Our present work
fits within this general context and is hence part of a large communal
push to develop scalable, efficient, robust data-driven algorithms for
analysis of large dynamical systems.

We proceed as follows. First, we briefly review some background
material on the basic theory of random Fourier features and the
N{\"y}strom method, and show how those ideas can be applied to EDMD. In the course of doing this, we also introduce some
thoughts on how to adaptively and efficiently add more features in the
random Fourier context in a way that makes use of previous
calculations. We then apply both the random Fourier and N{\"y}strom
techniques to two example EDMD problems, and show how these methods
provide an economical feature space for EDMD and hence an avenue for
analyzing data which has large state size and number of snapshots. We
conclude with an analysis of the computational runtime of the two
methods, and a comparison of the relative benefits and drawbacks of
each. We find that while the N{\"y}strom method is the more general and
accurate method, random Fourier features is the faster method for
larger problems.

\section{Random Fourier Features}
We first describe the random Fourier feature method and how it may
apply to EDMD. We begin with a brief review of the basic theory from
the literature, and conclude with a proposition on how to adaptively
add Fourier features in an EDMD problem efficiently.

\subsection{Basis Selection}
The main ideas we sketch here were first developed by Rahimi \&
Recht~\cite{Rahimi2007,Rahimi2008}. Coined ``Random Kitchen Sinks'',
the method developed therein sought to expand a function in terms of a
finite collection of random Fourier modes. The Fourier modes are
chosen in that algorithm according to a pre-specified frequency
distribution, which itself is derived from a pre-specified,
translation-invariant kernel function. Viewed another way, Random
Kitchen Sinks could be described as a method in which one first
chooses a translation-invariant kernel function (which computes inner
products in feature space), and then approximates that kernel function
by means of Monte Carlo sampling a Fourier basis. This is possible by
way of Bochner's theorem.

Bochner's theorem is the theoretical underpinning for the random
Fourier feature method. Consider the kernel function $k(\cdot,\cdot) :
\mathbb{R}^d \times \mathbb{R}^d \mapsto \mathbb{R}$. It states that
if that kernel function is positive semi-definite and
translation-invariant (i.e., $k(x,x') = k(x-x')$), then it may be
represented in a Fourier decomposition:
\begin{equation}
  \label{eq:bochner}
  k(x-x') = \int_z \psi_z(x) \overline{\psi_z(x')} \lambda(z) dz = \mathbb{E}_z \left[ \psi_z(x) \overline{\psi_z(x')} \right],
\end{equation}
where $x,x',z \in \mathbb{R}^d$, $\psi_z(x) = e^{i \langle z,x
  \rangle}$, and $\lambda(z)$ is a probability distribution associated
with the kernel. The random Fourier feature method seeks to
approximate this expectation with random sampling: a finite number $K$
of samples $\lbrace z_j \rbrace_{j=1}^K$ are drawn from the
distribution $\lambda(z)$, and the sample average is calcuated:
\begin{equation}
  \label{eq:rffsampleavg}
  k(x-x') \approx \frac{1}{K} \sum_{j=1}^K \psi_{z_j}(x) \overline{\psi_{z_j}(x')} \; .
\end{equation}
This approximation converges in expectation to the desired kernel; for
more information about the rate of convergence or other details of the
method, see~\cite{Rahimi2007,Rahimi2008,Le2013}.

The distribution $\lambda(z)$ may be found by taking the inverse
Fourier transform of $k(x)$, from which we can see that $\lambda(z)$
that depends on the form of $k(x)$. For example, given a Gaussian
kernel, the weights $\lambda(z)$ will also be normally-distributed in
the frequency domain (but with the inverse covariance). Other choices
of kernels (e.g., Laplacian, Cauchy) would lead to different
corresponding forms of $\lambda(z)$.

The application we propose in the context of EDMD is to simply use
this random Fourier basis as the EDMD features in
Eq.~\ref{eq:EDMDproblem}. 

\subsection{Adaptive Feature Addition}
\label{sec:AdaptiveRF}
One concern common to methods that approximate functions with a set of
basis functions is adaptivity: if we compute a function surrogate with
$K_0$ basis functions, and we wish to then add $K_{new}$ new basis
functions, what is the least computationally burdensome method for
doing this? In this section, we present a means for updating the
Koopman operator using a group of new basis functions in a way that
takes advantage of information from the previous computation and hence
provides a modest computational advantage to recomputing everything
from scratch. It should be noted that this method unfortunately does
not help us calculate the Koopman eigenspectrum any more quickly; it
is solely intended as a means to efficiently update the Koopman
operator.

The motivation for this partially stems from recent work on streaming
DMD~\cite{Hemati2014}. In that context, however, the problem was to
update the DMD Koopman approximation using sequentially obtained
snapshots of data. Our problem is to update the EDMD Koopman
approximation using a new group of basis functions (e.g., random
Fourier modes), and hence the tactics used are different.

Let $K_0$ be the initial number of basis functions in
Eq.~\ref{eq:GramKoopman}. Given that we wish to update the matrices
$\Psi_X$ and $\Psi_Y$ with new basis functions $\Psi_{X_{new}}$ and
$\Psi_{Y_{new}}$, our goal is to update the matrices $G$, $H$, and
$A$.

The first move we make is to simply note that since we have already
calculated the inner products of each of the old basis functions with
each other, we need not recalculate those terms in either $G$ or
$H$. The updated matrix $G$ can be written as:
\begin{equation}
\label{GramMatrix}
\begin{aligned}
G &= 
\left[
\begin{array}{c}
\Psi_{X_0}^T \\
\hline
\Psi_{X_{new}}^T
\end{array}
\right]
\left[
\begin{array}{c|c}
\Psi_{X_0} & \Psi_{X_{new}}
\end{array}
\right]
\\
&=
\left[
\begin{array}{c|c}
\Psi_{X_0}^T \Psi_{X_0} & \Psi_{X_0}^T \Psi_{X_{new}} \\
\hline
\Psi_{X_{new}}^T \Psi_{X_0} & \Psi_{X_{new}}^T \Psi_{X_{new}}
\end{array} 
\right] 
=
\left[
\begin{array}{c|c}
G_0 & G_1 \\
\hline
G_1^T & G_2
\end{array} 
\right] 
\end{aligned}
\end{equation}
The matrix $H$ has an analogous structure. Assuming we already have
knowledge of $G_0$ from the previous EDMD calculation, only need
calculate $G_1$ and $G_2$. The cost of calculating $G_1$ is
$\mathcal{O}(K_0 K_{new} M)$, and $G_2$ is $\mathcal{O}(K_{new}^2
M)$. Arguments for the cost of updating $H$ are exactly
parallel. Therefore, the cost of updating $G$ and $H$ is the
greater of either $\mathcal{O}(K_0 K_{new} M)$ or
$\mathcal{O}(K_{new}^2 M)$, depending on the relative sizes of $K_0$
or $K_{new}$, which denote the number of original and new basis
functions, respectively.

The next step involved is to update $G^\dagger$. It is shown
in~\cite{Rohde} that if $G$ is a nonnegative, symmetric matrix with
total rank equal to the sum of the ranks of $G_0$ and $G_2$, then the
pseudoinverse may be calculated in block form as:
\begin{equation}
\label{eq:PinvUpdate}
G^\dagger =
\left[
\begin{array}{c|c}
G_0^\dagger + G_0^\dagger G_1 Q^\dagger G_1^T G_0^\dagger & -G_0^\dagger G_1 Q^\dagger \\
\hline
-Q^\dagger G_1^T G_0^\dagger & Q^\dagger
\end{array} 
\right] 
\end{equation}
where:
\begin{equation}
\label{eq:PinvUpdate2}
Q = G_2 - G_1^T G_0^\dagger G_1
\end{equation}
Since $G$ is by definition a symmetric Gram matrix, a sufficient
condition for satisfying the assumptions needed for
Eq.~\ref{eq:PinvUpdate} is that the basis functions be linearly
independent. Furthermore, we already have knowledge of both $G_0$ and
$G_0^\dagger$, which alleviates the cost of calculating
Eq.~\ref{eq:PinvUpdate}. Therefore, the total cost of
Eq.~\ref{eq:PinvUpdate} is asymptotically dominated by the cost of the
matrix products involving $G_0^\dagger$ and $G_1$, which is
$\mathcal{O}(K_0^2 K_{new})$. The last step involves forming the
matrix product given in Eq.~\ref{eq:GramKoopman}, which requires
$\mathcal{O}((K_0 + K_{new})^3)$ time.

By comparison, explicit calculation from scratch of the updated
matrices $G$ and $H$ would require $\mathcal{O}((K_0+K_{new})^2 M)$
time; the pseudoinverse $G^\dagger$ would require
$\mathcal{O}((K_0+K_{new})^3)$ time; $G^\dagger H$ would require
$\mathcal{O}((K_0 + K_{new})^3)$ time. Therefore, the computational
savings obtained from using previous calculations would be largest
when the number of new basis functions added is small relative to the
number of original basis functions ($K_{new} \ll K_0$). Unfortunately,
the final asymptotic scaling of calculating the Koopman matrix
(Eq.~\ref{eq:GramKoopman}) is the same regardless of the method used,
but significant time could be saved in the calculation of $G$, $H$ and
$G^{\dagger}$.

\section{N{\"y}strom Approximation}
Thus far we have reviewed an approach to EDMD which could be described
as deductive: we exploited the a-priori knowledge that
translation-invariant kernels are well-approximated in the Fourier
basis to generate an efficient basis for EDMD. This connection between
a particular kernel and a corresponding basis is not always so clear,
however, particularly if the assumption of translation-invariance does
not hold. In light of this, one might wonder if there is a
complementary inductive approach, whereby one need only specify the
kernel function, and a basis for EDMD is then learned a-posteriori
from some training data. The N{\"y}strom
method~\cite{Nystrom1,Nystrom2} provides means for doing this.

The N{\"y}strom method is a data-driven approach toward approximating
the eigendecomposition of a kernel function. Here, we need not make
the assumption of translation-invariance of the kernel; we only assume
the kernel is symmetric and positive semi-definite. Under these
less restrictive conditions, Bochner's theorem no longer applies, but
the more general Mercer's theorem does. Mercer's theorem guarantees us
that any such kernel function $k(\cdot,\cdot) : \mathbb{R}^d \times
\mathbb{R}^d \mapsto \mathbb{R}$ may be expanded in terms of the
following infinite summation:
\begin{equation}
  k(x,y) = \sum_{j=1}^{\infty} \lambda_j \psi_j(x) \psi_j(y) \; ,
\end{equation}
where $\lbrace \lambda_j , \psi_j(x) \rbrace$ are
eigenvalues/functions of the Hilbert-Schmidt integral
operator that is associated with the kernel:
\begin{equation}
  \label{eq:KernelEigen}
  T_k \left( \psi_j(y) \right) = \int_x k(y,x) \psi_j(x) p(x) dx = \mu_j \psi_j(y) \; . 
\end{equation}
Here, $p(x)$ denotes the probability density function of the data
$x$. The N{\"y}strom method seeks to construct a Monte Carlo approximation
to the integral in Eq.~\ref{eq:KernelEigen} using a finite sample of data
$\lbrace x_1 \dots x_K \rbrace$ drawn from $p(x)$:
\begin{equation}
\frac{1}{K} \sum_j^K k(y,x_j) \psi_i(x_j) \approx \lambda_i \psi_i(y) \; .
\end{equation}
When this equation is written for each of the $K$ data points, we
produce the following matrix eigenproblem:
\begin{equation}
\label{eq:KernelMatrixEigen}
M_k U = U \Lambda
\end{equation}
where $M_k, U, \Lambda \in \mathbb{R}^{K \times K}$, $(M_k)_{i,j} =
k(x_i,x_j)$. Approximations of the kernel eigenvalues and
eigenfunctions at the $K$ data points are then:
\begin{equation}
\label{eq:NystromEigs}
\psi_i(x_j) \approx \sqrt{K} U_{j,i} \;\;\; , \;\;\; \lambda_i \approx \frac{\Lambda_{i,i}}{K}
\end{equation}
Furthermore, approximations at another data point $y$ may be interpolated:
\begin{equation}
\label{eq:NystromEigsInterp}
\psi_i(y) \approx \frac{\sqrt{K}}{\Lambda_{i,i}} \sum_j^K k(y,x_j) U_{j,i}
\end{equation}
To proceed for EDMD purposes, we must produce the feature matrices $\Psi_X$ and $\Psi_Y$
in Eq.~\ref{eq:EDMDKoopman}. At this point, we have a choice. One
option is to make use of the entire dataset and use
Eq.~\ref{eq:NystromEigsInterp} to interpolate $\Psi_X$ at the
remaining $\left(M-K\right)$ points and $\Psi_Y$ at all $M$ data points. A second,
cheaper option is to simply use the evaluation
Eq.~\ref{eq:NystromEigs} on the $K$ points that we already have as
$\Psi_X$, and only interpolate $\Psi_Y$ on those $K$ points using
Eq.~\ref{eq:NystromEigsInterp}. In what follows, the former method
will be referred to as the ``expensive'' N{\"y}strom variant; the latter
will be referred to as the ``cheap'' variant. Of course, one may also
elect to use an interpolation on some arbitrary subset of the data
between those two extremes. The choice is up to the user, although --
as we will see -- it can have predictable consequences in terms of the
speed/accuracy trade-off.

Notice that this N{\"y}strom method is more general than the random
Fourier feature method, as nowhere have we been required to assume
that the kernel function is translation-invariant. Although the
N{\"y}strom features are not theoretically exact as was the case in the
random Fourier feature method, they are derived from the data we
collect, which gives them a direct connection to the problem at hand
that may provide some compensating benefit.

\section{Numerical Examples}
We now present two numerical experiments that give some empirical
evidence demonstrating the efficacy of EDMD using either random
Fourier features or the N{\"y}strom method.

Before proceeding, the topic of kernel selection -- which is central
to both of these approaches -- should be discussed. Clearly, the
choice of a kernel can have a large effect on the efficiency and/or
accuracy of the results. The approach that we follow in this paper is
to use the data we have to estimate an acceptable kernel. We calculate
Euclidean distances between snapshots, and use the resulting
statistics to estimate the standard deviation of a Gaussian
distribution. It should be acknowledged up front that this is a bit of
a heuristic that may be unsatisfying to those who would want a more
rigorous approach. However, we were able to obtain satisfactory
results in the following examples with this. If a more rigorous
approach is desired, it may be possible to choose a kernel form (e.g.,
Gaussian), and then tune the shape parameter of that distribution with
cross-validation. One might even extend this cross-validation to
include multiple kernel forms (e.g., Laplacian, Cauchy) with various
shape parameters. It would be interesting to perform a follow-up study
to investigate whether an approach like this could yield improved
results; however, in order not to detract from the central narrative
of this paper, we do not pursue that discussion here any further.
  
\subsection{Fitzhugh-Nagumo Equations}
Here, following Williams~\cite{WilliamsKDMD_journal}, we examine the 1-D
Fitzhugh-Nagumo equations as a test case. This example is attractive
for several reasons. First, it provides a problem with a
moderately-sized state dimension, which makes it challenging for the
more traditional basis choices in EDMD (at least, without the
pre-application of some form of data compression, such as
POD). Second, the leading two Koopman modes can be deduced from the
system linearization and, hence, regular DMD can produce them for
comparison purposes.

The equations read thus:
\begin{equation}
\begin{aligned}
v_t &= v_{xx} + v - w - v^3 \\
w_t &= w_{xx} + \epsilon(v - c_1w - c_0)
\end{aligned}
\end{equation}
on the domain $x \in [0,20]$ and with the parameter values $c_0 =
-0.03$, $c_1 = 2.0$, $\delta = 4.0$, $\epsilon = 0.02$. The boundary
conditions used are Neumann, and we solve these equations using a
simple structured grid of 100 spatial points. Our goal is to extract
the Koopman modes for $v$, and so our state dimension $d = 100$. We
use 2500 snapshot pairs, taken at temporal intervals of $\Delta t =
1$. The initial conditions used are the standing wave fronts in $v$
and $w$, and the system is perturbed every $25\Delta t$ by a simple
Gaussian process.

\begin{figure}[!htb]
\centering
\begin{subfigure}[t]{0.6\textwidth}
\includegraphics[width=1\textwidth]{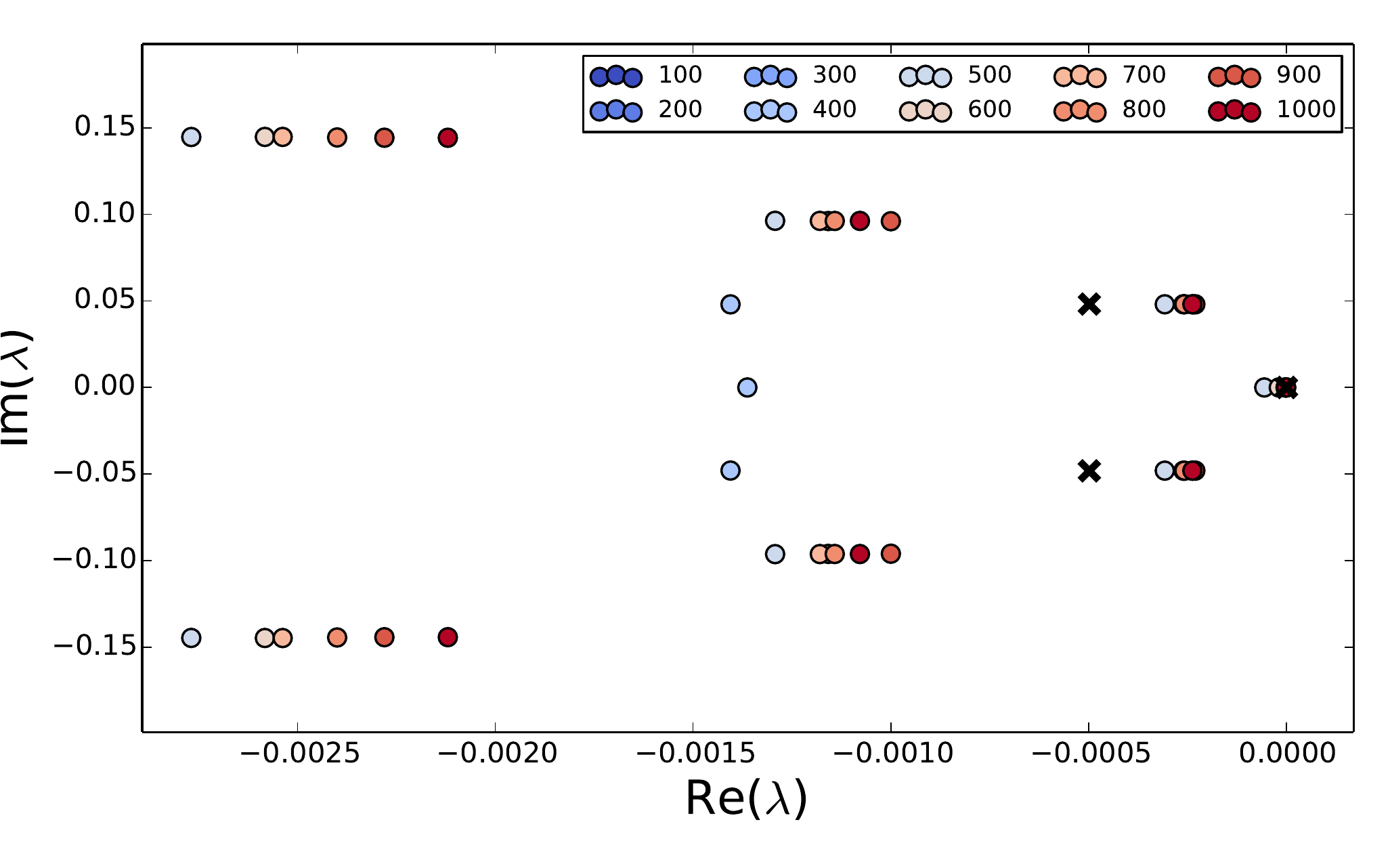}
\caption{Eigenvalues (DMD eigenvalues displayed as black $\times$ markers).}
\label{fig:FNeigA}
\end{subfigure}\\
\begin{subfigure}[t]{0.6\textwidth}
\includegraphics[trim={0 2.5cm 0 0},clip,width=1\textwidth]{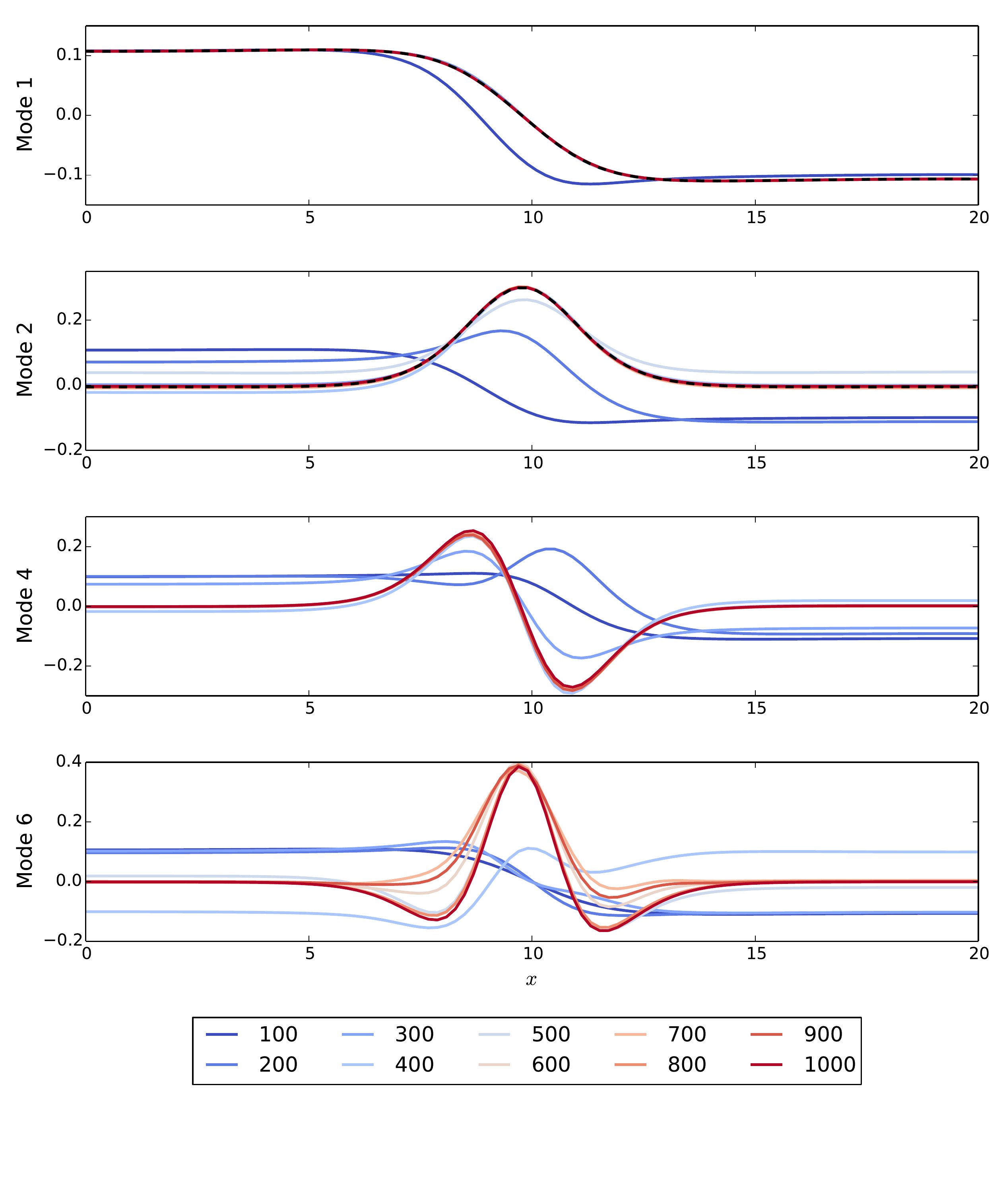}
\caption{Modes (DMD modes displayed as dashed black lines).}
\label{fig:FNeigB}
\end{subfigure}
\caption{Fitzhugh-Nagumo eigenvalues and modes, using DMD and random Fourier EDMD with varying numbers of basis functions.}
\label{fig:FNeig}
\end{figure}

\begin{figure}[!htb]
\centering
\begin{subfigure}[t]{0.6\textwidth}
\includegraphics[width=1\textwidth]{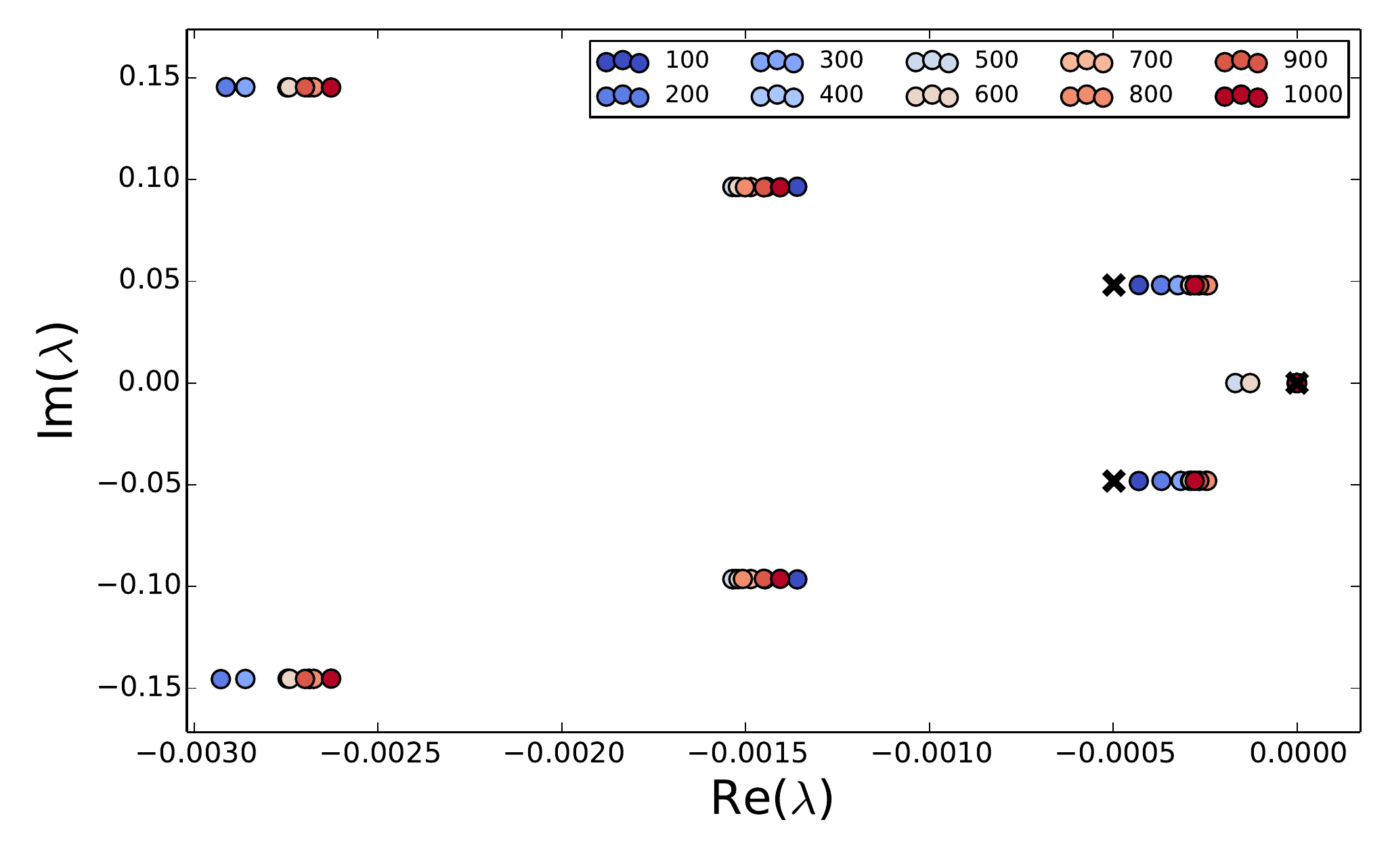}
\caption{Eigenvalues (DMD eigenvalues displayed as black $\times$ markers).}
\label{fig:FNeig2A}
\end{subfigure}\\
\begin{subfigure}[t]{0.6\textwidth}
\includegraphics[trim={0 2.5cm 0 0},clip,width=1\textwidth]{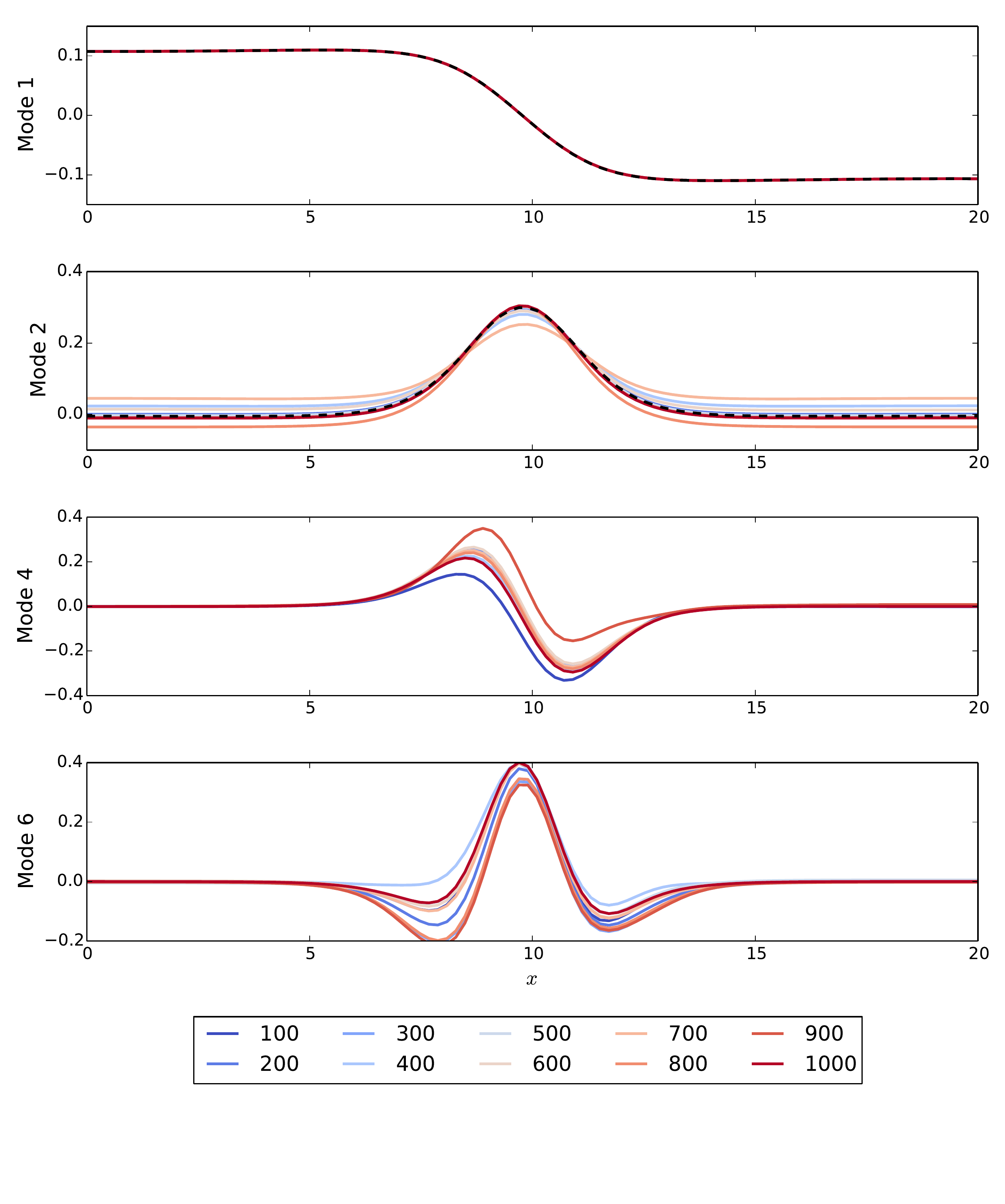}
\caption{Modes (DMD modes displayed as dashed black lines).}
\label{fig:FNeig2B}
\end{subfigure}
\caption{Fitzhugh-Nagumo eigenvalues and modes, using DMD and the ``expensive'' N{\"y}strom method EDMD with varying numbers of basis functions.}
\label{fig:NystromFNeig}
\end{figure}

Fig.~\ref{fig:FNeig} shows computations of the leading Koopman
eigenvalues and modes using a random Fourier basis with Monte Carlo
sample sizes ranging from 100 to 1000. Fig.~\ref{fig:NystromFNeig}
shows the equivalent computations using the expensive N{\"y}strom method
with random snapshot sample sizes ranging from 100 to 1000. The random
Fourier basis frequencies are normally-distributed with a standard deviation of $\sigma =
4\pi$ (corresponding to a normally-distributed kernel function with
the inverse standard deviation); the kernel function used in the
N{\"y}strom method is a Gaussian with $\sigma = 1$. In both cases, we
determine the shape parameter $\sigma$ empirically from the data set:
we calculate distances between data snapshots and estimate an average
distance between them.

In the random Fourier feature results, the two leading (linear)
Koopman modes converge relatively quickly and strongly to the correct
answer (as judged by regular DMD), requiring around 500-600 basis
functions to accurately represent. The corresponding Koopman
eigenvalues converge quickly as well. The next two (nonlinear) Koopman
modes present, predictably, more variance in shape with number of
basis functions, and require more basis functions (around 900) to
converge. Similarly, there is more variation in the calculated Koopman
eigenvalues that correspond to those two nonlinear Koopman modes.

It appears Koopman modes computed with the expensive N{\"y}strom method
may be converging sooner than those computed with random Fourier
features. The most noticeable difference demonstrating this is the
shape of the modes using lower numbers of snapshots. In the Fourier
approach, around 500-600 basis functions were required before the
modes began to qualitatively resemble the correct answers; in the
N{\"y}strom approach, qualitative correctness is achieved almost instantly
at 100 snapshots.

It should be clear, however, that both methods provide fairly good
approximations of the leading Koopman modes/eigenvalues in this
problem, with only 1000 random samples or less. This is a reasonable
number of basis functions for a state dimension of $d = 100$, compared
to the more traditional basis choices.

\subsection{Experimental Cylinder Flow Data}
Our final demonstration involves a data set in which both the state
dimension and the number of snapshots are both quite large. We analyze
PIV data of 2D cylinder flow~\cite{TuShangPIV_journal}, consisting of a
$80 \times 135$ Cartesian grid of spatial points ($d = 10,800$)
sampled at $M = 8,000$ snapshots in time. This data set has been
analyzed previously using DMD~\cite{TuShangPIV_journal},
KDMD~\cite{WilliamsKDMD_journal}, and streaming DMD~\cite{Hemati2014}; the
interested reader is encouraged to consult those references for more
information on this data set.

As both $d$ and $M$ are approximately $10^4$, both regular DMD and
EDMD (with traditional basis choices) are either impractical or
impossible for the full data set: DMD can be performed on this data
set (see~\cite{TuShangPIV_journal}), but that requires several cores running
in parallel and is extremely computationally expensive; EDMD based on
tensor polynomials or RBFs clearly suffers from the large state
dimensionality. KDMD necessitates the SVD computation of a $M \times
M$ matrix ($\mathcal{O}(M^3)$ computational scaling) and hence is
expensive as well; because of this expense, it is only performed on a
subset of the full data in~\cite{WilliamsKDMD_journal}. Of course, it may be
argued that a pre-processing data compression step (e.g., POD) could
be applied to make the problem more computationally tractable by
reducing the effective state dimension. While that could be a feasible
approach, we are interested in addressing situations where the state
dimension (and number of snapshots) is large (which might be the case
even after any data pre-processing, depending on the
problem). Therefore, as in the previous example, we do not apply any
such dimension-reduction pre-processing.

This data set is well-suited for analysis using either random Fourier
EDMD or the N{\"y}strom method. As in the previous example, we select
the Fourier modal frequencies from a normal distribution, which
effectively approximates a Gaussian kernel with the inverse
variance. The Fourier modal distribution variance is, as before,
determined empirically from the data set (we use $\sigma = 1/400$). In
the N{\"y}strom approach, we use a Gaussian kernel function with the
inverse shape parameter used in the random Fourier method (i.e., the
N{\"y}strom kernel is a Gaussian with $\sigma = 400$). We test both
the ``expensive'' version of the N{\"y}strom method, in which we
interpolate the feature matrices $\Psi_X$ and $\Psi_Y$ on all $M$ data
points, and the ``cheap'' version, in which we only interpolate
$\Psi_Y$ on the $K$ random sample data points.

\begin{figure}[!htb]
\centering
\begin{subfigure}[t]{1\textwidth}
\includegraphics[trim={1.25in 1in 1in 0.75in},clip,width=1\textwidth]{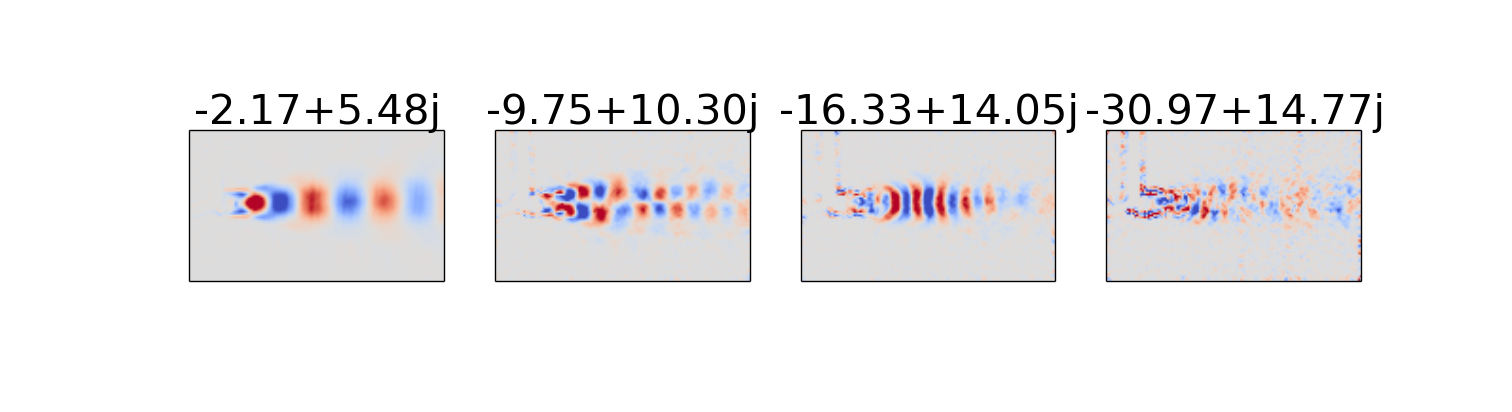}
\caption{$K=100$}
\label{fig:PIVa}
\end{subfigure}\\
\begin{subfigure}[t]{1\textwidth}
\includegraphics[trim={1.25in 1in 1in 0.75in},clip,width=1\textwidth]{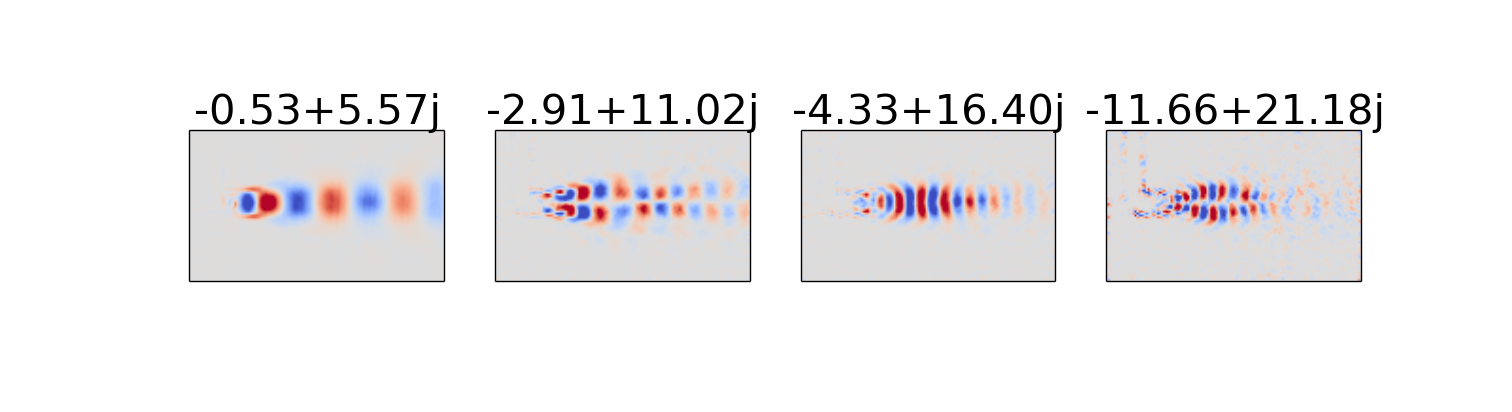}
\caption{$K=500$}
\label{fig:PIVb}
\end{subfigure}\\
\begin{subfigure}[t]{1\textwidth}
\includegraphics[trim={1.25in 1in 1in 0.75in},clip,width=1\textwidth]{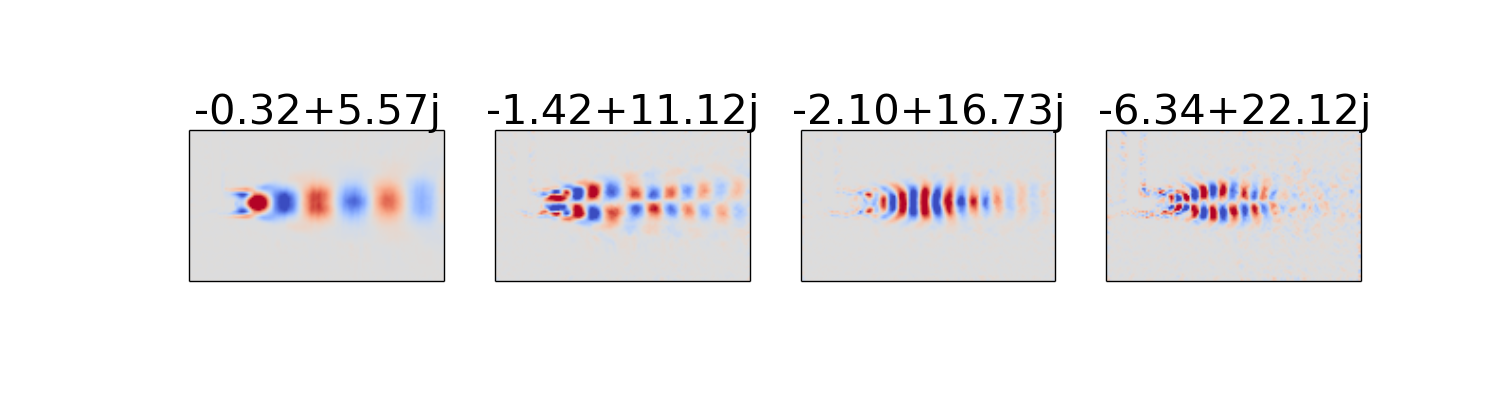}
\caption{$K=1000$}
\label{fig:PIVc}
\end{subfigure}\\
\begin{subfigure}[t]{1\textwidth}
\includegraphics[trim={1.25in 1in 1in 0.75in},clip,width=1\textwidth]{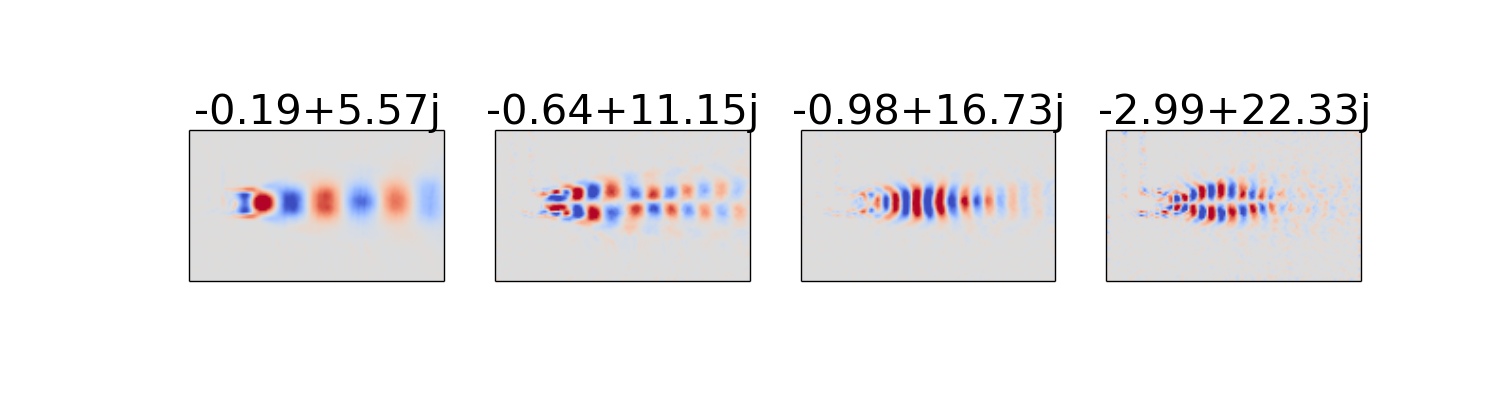}
\caption{$K=2000$}
\label{fig:PIVd}
\end{subfigure}
\caption{Leading Koopman mode and eigenvalue approximations for the cylinder PIV data, using random Fourier EDMD with varying numbers of basis functions ($K$).}
\label{fig:PIV}
\end{figure}

\begin{figure}[!htb]
\centering
\begin{subfigure}[t]{1\textwidth}
\includegraphics[trim={1.25in 1in 1in 0.75in},clip,width=1\textwidth]{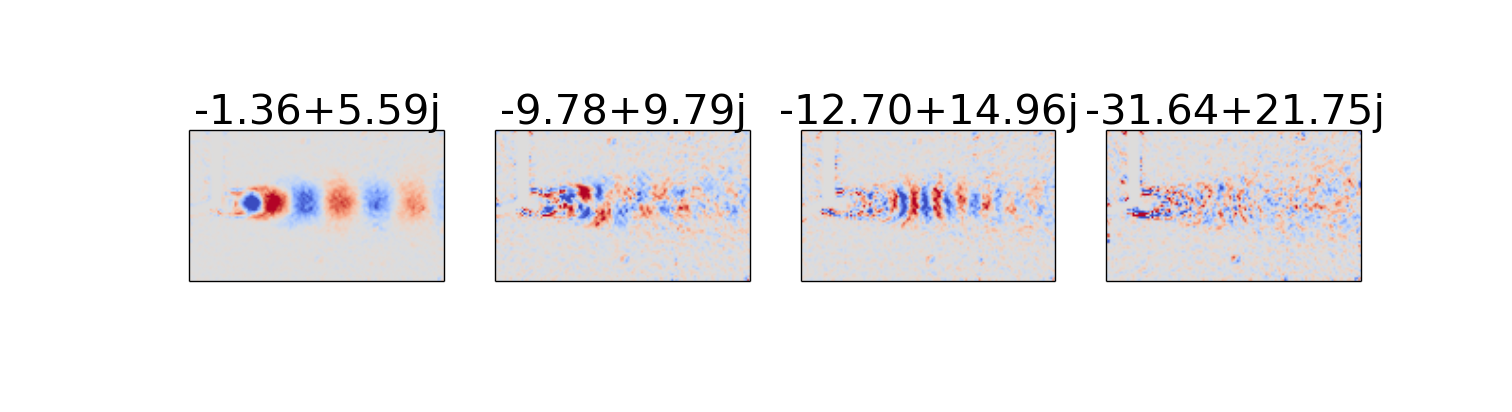}
\caption{$K=100$}
\label{fig:PIV2a}
\end{subfigure}\\
\begin{subfigure}[t]{1\textwidth}
\includegraphics[trim={1.25in 1in 1in 0.75in},clip,width=1\textwidth]{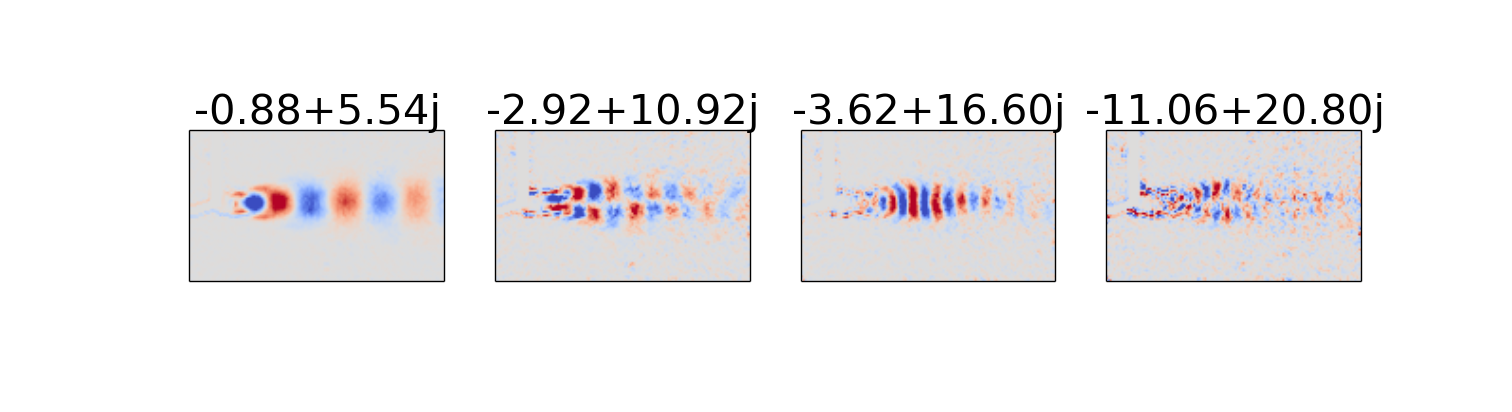}
\caption{$K=500$}
\label{fig:PIV2b}
\end{subfigure}\\
\begin{subfigure}[t]{1\textwidth}
\includegraphics[trim={1.25in 1in 1in 0.75in},clip,width=1\textwidth]{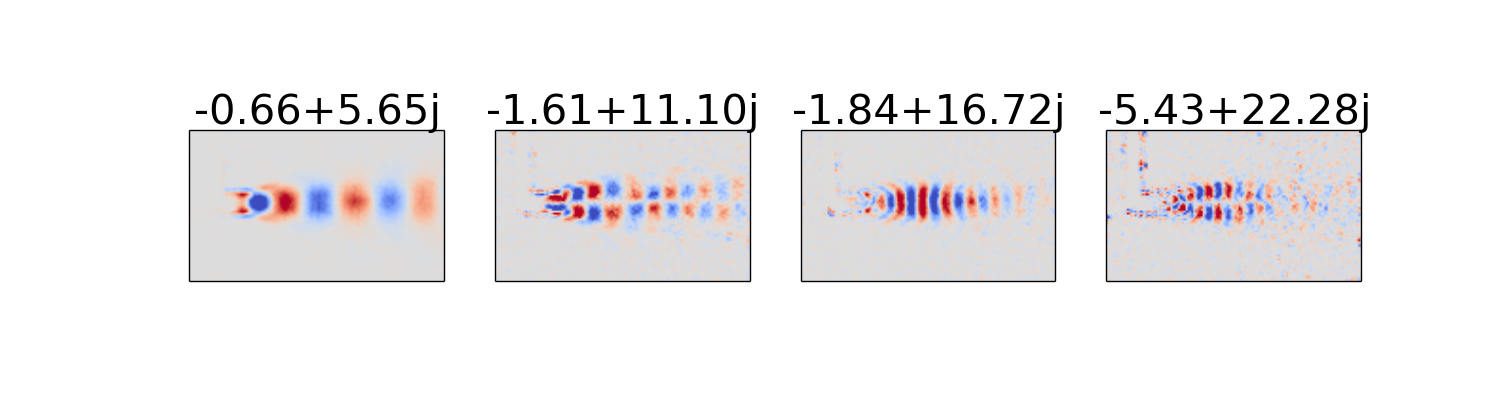}
\caption{$K=1000$}
\label{fig:PIV2c}
\end{subfigure}\\
\begin{subfigure}[t]{1\textwidth}
\includegraphics[trim={1.25in 1in 1in 0.75in},clip,width=1\textwidth]{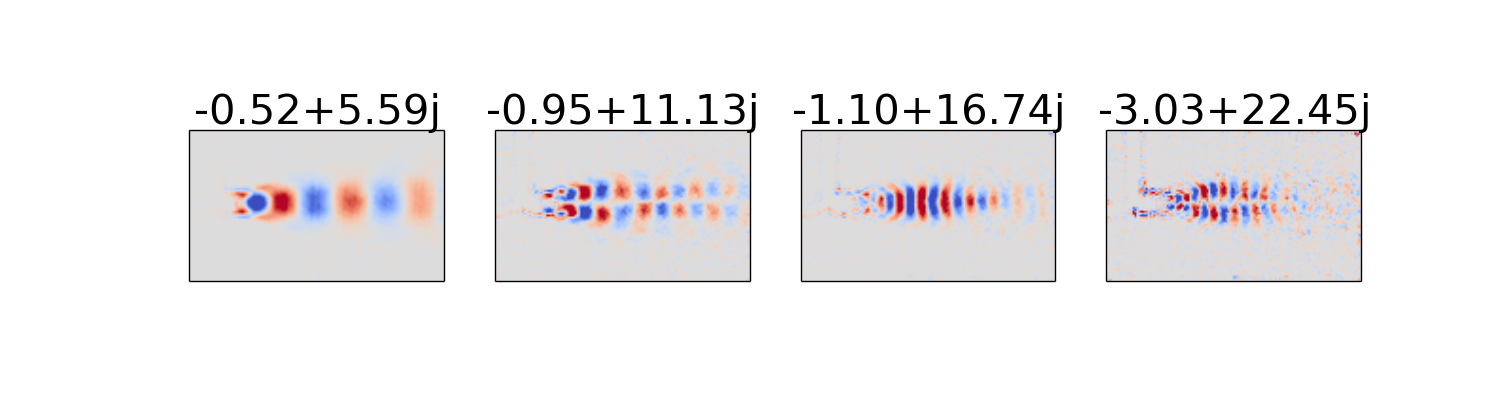}
\caption{$K=2000$}
\label{fig:PIV2d}
\end{subfigure}
\caption{Leading Koopman mode and eigenvalue approximations for the cylinder PIV data, using the ``cheap'' N{\"y}strom method EDMD with varying numbers of basis functions ($K$).}
\label{fig:PIV2}
\end{figure}

\begin{figure}[!htb]
\centering
\begin{subfigure}[t]{1\textwidth}
\includegraphics[trim={1.25in 1in 1in 0.75in},clip,width=1\textwidth]{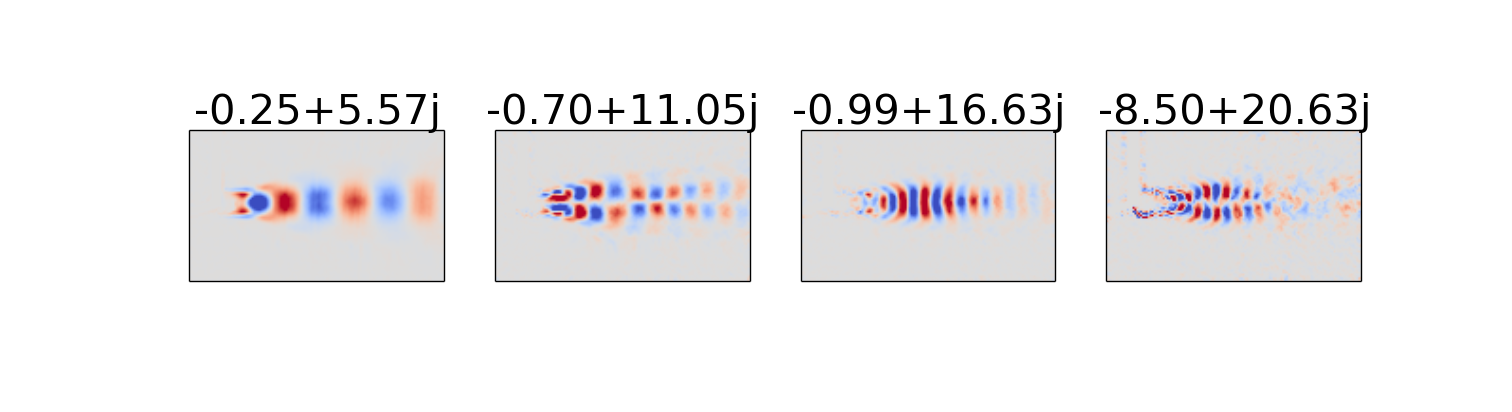}
\caption{$K=100$}
\label{fig:PIV3a}
\end{subfigure}\\
\begin{subfigure}[t]{1\textwidth}
\includegraphics[trim={1.25in 1in 1in 0.75in},clip,width=1\textwidth]{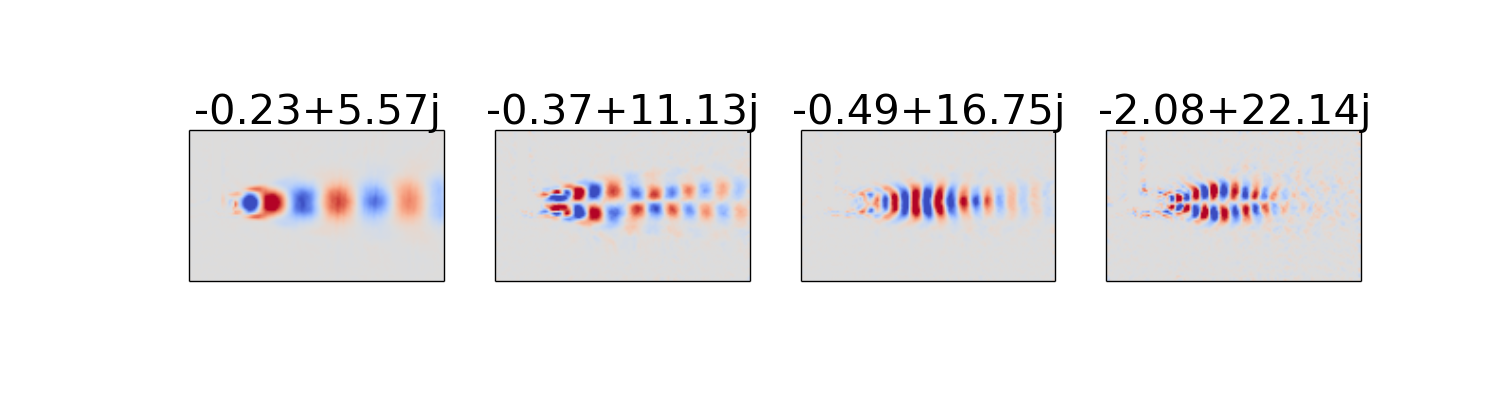}
\caption{$K=500$}
\label{fig:PIV3b}
\end{subfigure}\\
\begin{subfigure}[t]{1\textwidth}
\includegraphics[trim={1.25in 1in 1in 0.75in},clip,width=1\textwidth]{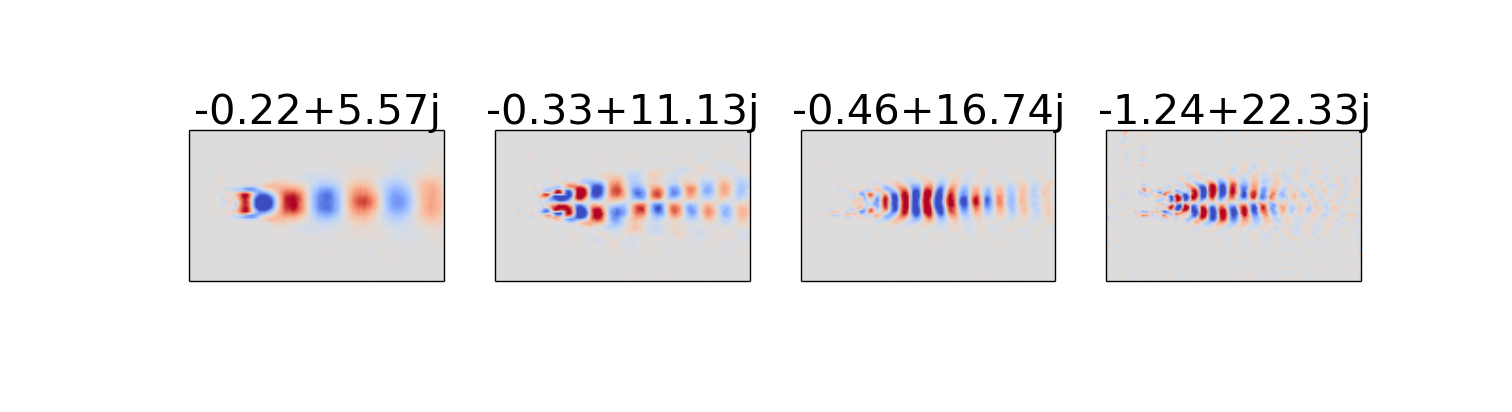}
\caption{$K=1000$}
\label{fig:PIV3c}
\end{subfigure}\\
\begin{subfigure}[t]{1\textwidth}
\includegraphics[trim={1.25in 1in 1in 0.75in},clip,width=1\textwidth]{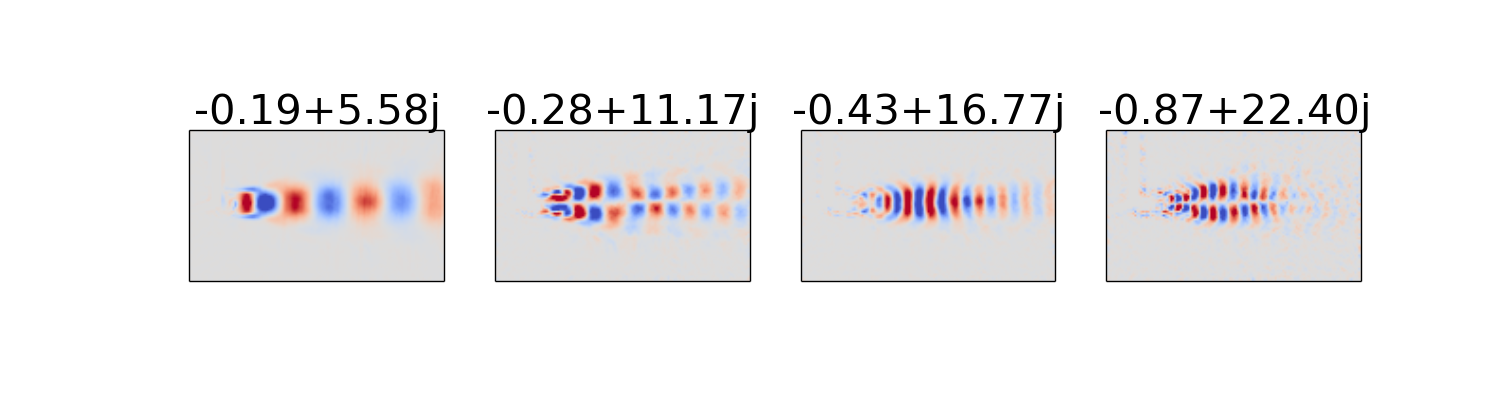}
\caption{$K=2000$}
\label{fig:PIV3d}
\end{subfigure}
\caption{Leading Koopman mode and eigenvalue approximations for the cylinder PIV data, using the ``expensive'' N{\"y}strom method EDMD with varying numbers of basis functions ($K$).}
\label{fig:PIV3}
\end{figure}

Fig.~\ref{fig:PIV} displays the leading Koopman modes and eigenvalues,
computed with the random Fourier method, for $K=100,\;200,\;500$, and
$1000$. Figs.~\ref{fig:PIV2} and~\ref{fig:PIV3} show equivalent
computations using the cheap and expensive N{\"y}strom variants
(respectively) with identical numbers $K$ of features. Using the
Fourier method, we see good convergence of the first four Koopman
modes when using around $500$ to $1000$ Fourier modes; this is
remarkable given both the large state and snapshot sizes of the data
set. Certainly -- as before -- the higher order modes require more
basis functions to converge (e.g., the first mode converges well using
$K < 100$, while the fourth Koopman mode requires 500 to 1000 basis
functions), but this asymetry in modal resolution is found in
EDMD/KDMD as well and so does not represent a weakness unique to
Fourier EDMD. We should also note that the Fourier EDMD method (as
well as the N{\"y}strom method) retains an advantage of EDMD/KDMD in that
the leading Koopman modes may be determined by the proximity of their
eigenvalues to the imaginary axis (as opposed to the ``cloud'' of
eigenvalues produced in DMD in which the dominant modes must be
determined using energy-based or sparsity-promoting methods).

In comparison, it appears yet again that convergence is relatively
faster using the expensive N{\"y}strom variant, with respect to both the
eigenvalues and eigenvectors. In particular, the fourth Koopman mode
is qualitatively well approximated in the expensive N{\"y}strom method
using only 100 empirical eigenfunction features, while 100 random
Fourier features is not sufficient for this purpose. However, accuracy
using the cheap N{\"y}strom variant is more or less equal to that using
random Fourier features.

In all methods, it does appear that convergence is slowest for the
real parts of the Koopman eigenvalues -- the four modes should all
have eigenvalues lying on the imaginary axis. While the imaginary
components are in good agreement with the reported values in previous
studies, the real components are noticeably larger in magnitude than
they should be (although this magnitude clearly decreases with
increasing $K$). One consequence of this is a ``bowing'' of the
eigenvalues, which is a common observation in EDMD/KDMD. Regardless of
this, we clearly see the utility of the two approaches on display in
this example: we are able to extract good approximations of the
leading components of the Koopman spectrum using only around 1000 EDMD
features, which is in stark contrast to the time and resources
required to do the same using DMD, KDMD, or EDMD (with traditional
basis choices).

\subsection{Computational Scaling}
Computationally speaking, there are some differences between the
random Fourier, cheap N{\"y}strom, and expensive N{\"y}strom methods. The
first difference involved is in how the basis elements are
computed. Basis computation in the N{\"y}strom method occurs in two
distinct steps. The first is to compute the empirical eigenfunctions
at the $K$ random data points. This involves computing the kernel
matrix for the $K$ random sample points ($\mathcal{O}(\frac{1}{2}K^2
d)$) and then computing the eigendecomposition of that $K \times K$
kernel matrix ($\mathcal{O}(K^3)$). The second is to interpolate those
eigenfunctions. In the expensive N{\"y}strom variant, this interpolation
occurs for $\Psi_X$ and $\Psi_Y$ at all $M$ data points
(Eq.~\ref{eq:NystromEigsInterp}), which involves computing two
$M \times K$ kernel matrices ($\mathcal{O}(M K d)$) and multiplying
each of them with a $K \times K$ matrix ($\mathcal{O}(2 M
K^2)$). Thus, the total runtime for the expensive N{\"y}strom basis
computation scales asymptotically as $\mathcal{O}(KMd + K^2(d+M) +
K^3)$. In the cheap N{\"y}strom variant, eigenfunction interpolation only
occurs for $\Psi_Y$ at the $K$ random data points, which involves
computing a $K \times K$ kernel matrix ($\mathcal{O}(\frac{1}{2} K^2
d)$) and multiplying it with a $K \times K$ matrix
($\mathcal{O}(K^3)$). Thus, the asymptotic scaling for the cheap
N{\"y}strom variant is $\mathcal{O}(K^2d + K^3)$. In random Fourier EDMD,
the random Fourier basis must simply be calculated at all data points,
which is $\mathcal{O}(KMd)$. Thus, we see that the random Fourier
method will always be faster than the expensive N{\"y}strom method for
basis generation; how the random Fourier method compares to the cheap
N{\"y}strom method will generally depend on the parameter values.

The remaining two steps involved in EDMD are the calculation of the
Koopman matrix and its eigendecomposition (for details,
see~\cite{WilliamsEDMD_journal}). The Koopman operator calculation scales as
$\mathcal{O}(K^2M + K^3)$. The eigendecomposition calculation scales
as $\mathcal{O}(KMd + K^2d + K^3)$ for the random Fourier and
expensive N{\"y}strom methods and $\mathcal{O}(K^2d + K^3)$ for the cheap
N{\"y}strom method. These asymptotic scalings are summarized in
Table~\ref{table:Scalings}. Note that in this table, ``CN'' denotes
the cheap N{\"y}strom method, ``RF'' the random Fourier method, and ``EN''
the expensive N{\"y}strom method.

\begin{table}[ht]
\caption{Asymptotic Computational Scalings}
\centering 
\begin{tabular}{c c c c c c c}
\hline
\hline
 & \vline & Basis & Koopman & Eigenspectrum \\ [0.5ex]
\hline
CN    & \vline & $K^2d + K^3$            & $K^3$         & $K^2d + K^3$          \\
RF    & \vline & $KMd$                  & $K^2M + K^3$   & $KMd + K^2d + K^3$    \\
EN    & \vline & $KMd + K^2(d+M) + K^3$  & $K^2M + K^3$  & $KMd + K^2d + K^3$     \\[1ex]
\hline
\end{tabular}
\label{table:Scalings}
\end{table}

\begin{figure}[!ht]
\centering
\scalebox{0.6}{\input{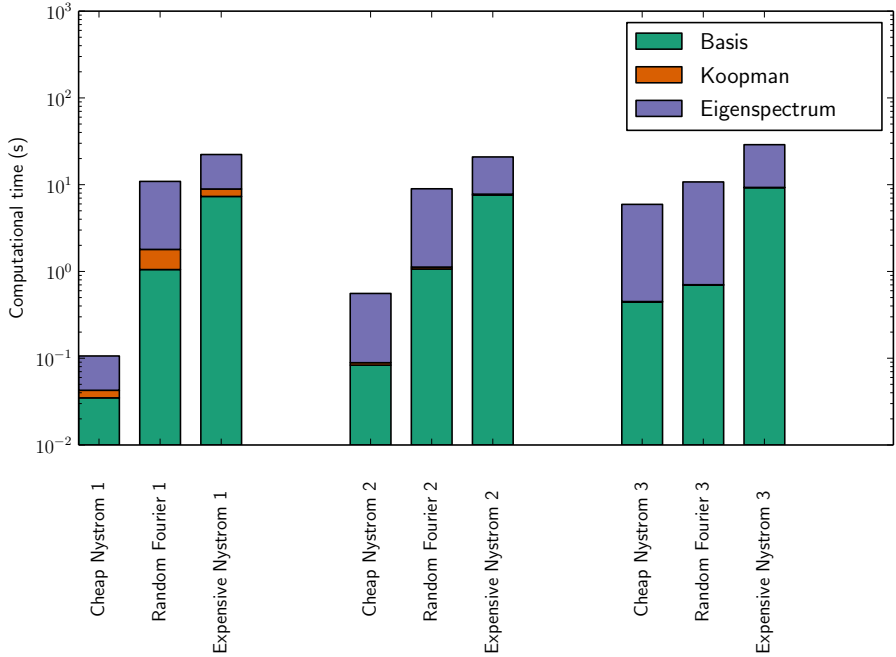}}
\caption{Computational runtimes. Parameters used in the five tests are: $d = \lbrace10^3, 10^4, 10^5\rbrace$, $M = \lbrace10^5, 10^4, 10^3\rbrace$, $K = \lbrace10^2, 10^2, 10^2\rbrace$, respectively. Bar color divisions are linear percentages of the total respective computational runtimes (given as the total bar heights).}
\label{fig:RunTimes}
\end{figure}

To illustrate concretely by way of examples, Fig.~\ref{fig:RunTimes}
presents runtimes for the random Fourier and N{\"y}strom methods using
problems with increasing values of $K$, $M$, and $d$, broken down by
each of the steps in the process (basis computation, Koopman matrix
calculation, and Koopman eigendecomposition). The three tests
displayed in that figure were chosen to investigate runtimes in the
regime where $K \ll d,M$, since that is the regime in which these
methods are intended to be used. The total runtimes of both the random
Fourier and expensive N{\"y}strom methods remain relatively constant for
each of the tests, while that of the cheap N{\"y}strom method increases
linearly with increasing $d$. This is the expected behavior on the
basis of Table~\ref{table:Scalings}: the product $KMd$ is the dominant
term in the scaling laws for the random Fourier and expensive N{\"y}strom
methods, and it is constant throughout each of the tests. Conversely,
the cheap N{\"y}strom method does not scale with $KMd$; the linear upward
trend visible is a result of the $K^2d$ term in both the basis and
eigendecomposition calculations. We also see clear confirmation that
the expensive N{\"y}strom method generally has the longest total runtime
of all the methods, a consequence of the high cost of basis
computation. While the cheap N{\"y}strom method is less expensive than the
random Fourier method in these three examples, it is important to note
that we cannot in general conclude that it will always be faster. We
can make the general observation that the cheap N{\"y}strom method will be
the fastest when $d \ll M$ (i.e., lower state dimension, higher number
of snapshots), since the random Fourier and expensive N{\"y}strom methods
scale linearly with $M$ but the cheap N{\"y}strom method is independent of
$M$ (a nice feature).

Incorporating this information yields a complex set of relative
benefits and drawbacks when comparing the methods. Our numerical
experiments seem to indicate that the expensive N{\"y}strom method is more
accurate than the other methods, given the same values of
$K,M,d$. Both of the N{\"y}strom methods are more flexible, in the sense
that they can approximate a wider class of kernels, and they have a
closer relationship to the data. Additionally, the cheap N{\"y}strom
method does not scale with $M$, since it only uses a random
subsampling of the full data set to generate the empirical
basis. However, the random Fourier method uses eigenfunctions that are
more mathematically precise for translation-invariant kernels than the
data-approximated ones generated by N{\"y}strom. Additionally, we have a
method for efficiently adding random Fourier features
(\S~\ref{sec:AdaptiveRF}), which is another speed-based attractive
feature.

Lastly, some comments on how these methods compare to ``regular''
EDMD/KDMD are in order. It is important here to recall the overall
objective of this paper: our goal was not to improve the asymptotic
computational scaling of EDMD with respect to the parameters $K,M,d$;
rather, it was to develop and use highly efficient EDMD basis
functions that effectively make the value of $K$ much less than what
it otherwise would be using traditional basis choices for problems
with large $d,M$. Hence, although the computational scaling for
regular EDMD is essentially the same as it is for the N{\"y}strom and
random Fourier methods, both the N{\"y}strom and random Fourier methods
are able to handle problems that would require an infeasibly large
number of basis functions if approached with regular EDMD. Similar
arguments exist regarding KDMD, as its runtime is dominated by $M^3$.

\section{Conclusions}
The objective of this work was to make progress toward reducing the
computational demands of EDMD/KDMD for data sets with large state and
snapshot sizes. To do this, we imported ideas from the theory of
random approximations of kernel functions, and used these ideas to
give us an economical set of EDMD features. We presented two numerical
examples which demonstrate how we can obtain good approximations of
the Koopman modes/eigenvalues for problems with large state and
snapshot sizes using a very reasonable number of basis functions.

We observed a complex set of trade-offs involved in comparing the
random Fourier and N{\"y}strom methods. We saw that the expensive N{\"y}strom
method generally provides higher accuracy, but at an increased
computational cost, for the same values of $K$, $M$, and $d$. The
cheap N{\"y}strom method computational runtime does not scale with
$M$. The N{\"y}strom methods are more flexible in the sense that they can
represent a wider range of kernels, while the random Fourier method is
more flexible in the sense that features can be added adaptively.

One avenue that may be of future research interest involves speeding
up the random Fourier basis computation. It has already been shown
in~\cite{Le2013} that Hadamard transforms may be used to obtain fast
approximations to the random Fourier basis computation; however, the
method as it stands computes a number of basis functions that is equal
to or greater than the state dimension, which might cancel its
economical gains when $K \ll d$ (as is the case in many problems of
interest).

\section{Acknowledgments}
Los Alamos Report LA-UR-17-29880. Funded by the Department of Energy
at Los Alamos National Laboratory under contract DE-AC52-06NA25396. The authors also wish to acknowledge Scott Dawson for helpful
technical feedback and Mark Lohry for helpful advice on visually displaying data.

\clearpage

\bibliographystyle{siamplain}
\bibliography{RandomFourierEDMD}

\begin{thebibliography}{10}

\bibitem{bagheri}
{\sc S.~Bagheri}, {\em Effects of weak noise on oscillating flows: linking
  quality factor, floquet modes, and koopman spectrum}, Physics of Fluids, 26
  (2014), p.~094104.

\bibitem{bollt}
{\sc E.~M. Bollt and N.~Santitissadeekorn}, {\em Applied and Computational
  Measurable Dynamics}, vol.~18, SIAM, 2013.

\bibitem{BruntonCompressiveDMD}
{\sc S.~L. Brunton, J.~L. Proctor, J.~H. Tu, and J.~N. Kutz}, {\em {Compressive
  sampling and dynamic mode decomposition}}, Journal of Computational Dynamics,
  2 (2018), pp.~165--191.

\bibitem{budisic}
{\sc M.~Budi\v{s}i\'{c}, R.~Mohr, and I.~Mezi\'{c}}, {\em Applied koopmanism},
  Chaos: An Interdisciplinary Journal of Nonlinear Science, 22 (2012),
  p.~047510.

\bibitem{chen}
{\sc K.~K. Chen, J.~H. Tu, and C.~W. Rowley}, {\em Variants of dynamic mode
  decomposition: boundary condition, koopman, and fourier analyses}, Journal of
  Nonlinear Science, 22 (2012), pp.~887--915.

\bibitem{coifman}
{\sc R.~R. Coifman and S.~Lafon}, {\em Diffusion maps}, Applied and
  Computational Harmonic Analysis, 21 (2006), pp.~5--30.

\bibitem{dahl}
{\sc G.~E. Dahl, D.~Yu, L.~Deng, and A.~Acero}, {\em Context-dependent
  pre-trained deep neural networks for large-vocabulary speech recognition},
  IEEE Transactions on audio, speech, and language processing, 20 (2012),
  pp.~30--42.

\bibitem{dawson}
{\sc S.~T. Dawson, M.~S. Hemati, M.~O. Williams, and C.~W. Rowley}, {\em
  Characterizing and correcting for the effect of sensor noise in the dynamic
  mode decomposition}, Experiments in Fluids, 57 (2016).

\bibitem{ErichsonCompressedDMD}
{\sc N.~B. Erichson, S.~L. Brunton, and J.~N. Kutz}, {\em {Compressed dynamic
  mode decomposition for background modeling}}, Journal of Real-Time Image
  Processing,  (2016), pp.~1--14.

\bibitem{ErichsonRandomizedDMD2017}
{\sc N.~B. Erichson, S.~L. Brunton, and J.~N. Kutz}, {\em {Randomized dynamic
  mode decomposition}}.
\newblock arXiv:1702.02912, 2017.

\bibitem{ErichsonRandomizedDMD2016}
{\sc N.~B. Erichson and C.~Donovan}, {\em {Randomized low-rank dynamic mode
  decomposition for motion detection}}, Computer Vision and Image
  Understanding, 146 (2016), pp.~40--50.

\bibitem{froyland}
{\sc G.~Froyland, G.~A. Gottwald, and A.~Hammerlindl}, {\em A computational
  method to extract macroscopic variables and their dynamics in multiscale
  systems}, SIAM Journal on Applied Dynamical Systems, 13 (2014),
  pp.~1816--1846.

\bibitem{girshick}
{\sc R.~Girshick, J.~Donahue, T.~Darrell, and J.~Malik}, {\em Rich feature
  hierarchies for accurate object detection and semantic segmentation}, in
  Proceedings of the IEEE Conference on Computer Vision and Pattern
  Recognition, 2014, pp.~580--587.

\bibitem{gueniat2015}
{\sc F.~Gu{\'e}niat, L.~Mathelin, and L.~R. Pastur}, {\em A dynamic mode
  decomposition approach for large and arbitrarily sampled systems}, Physics of
  Fluids, 27 (2015), p.~025113.

\bibitem{hemati2017}
{\sc M.~S. Hemati, C.~W. Rowley, E.~A. Deem, and L.~N. Cattafesta}, {\em
  De-biasing the dynamic mode decomposition for applied koopman spectral
  analysis of noisy datasets}, Theoretical and Computational Fluid Dynamics, 31
  (2017), pp.~349--368.

\bibitem{Hemati2014}
{\sc M.~S. Hemati, M.~O. Williams, and C.~W. Rowley}, {\em {Dynamic mode
  decomposition for large and streaming datasets}}, Physics of Fluids, 26
  (2014).

\bibitem{hinton}
{\sc G.~Hinton, L.~Deng, D.~Yu, G.~E. Dahl, A.~Mohamed, N.~Jaitly, A.~Senior,
  V.~Vanhoucke, P.~Nguyen, T.~N. Sainath, and B.~Kingsbury}, {\em Deep neural
  networks for acoustic modeling in speech recognition: The shared views of
  four research groups}, IEEE Signal Processing Magazine, 29 (2012),
  pp.~82--97.

\bibitem{jovanovic}
{\sc M.~R. Jovanovi{\'c}, P.~J. Schmid, and J.~W. Nichols}, {\em
  Sparsity-promoting dynamic mode decomposition}, Physics of Fluids, 26 (2014),
  p.~024103.

\bibitem{koopman}
{\sc B.~O. Koopman}, {\em Hamiltonian systems and transformation in hilbert
  space}, Proceedings of the National Academy of Sciences, 17 (1931),
  pp.~315--318.

\bibitem{krizhevsky}
{\sc A.~Krizhevsky, I.~Sutskever, and G.~E. Hinton}, {\em Imagenet
  classification with deep convolutional neural networks}, Advances in Neural
  Information Processing Systems,  (2012), pp.~1097--1105.

\bibitem{Le2013}
{\sc Q.~Le, T.~Sarlos, and A.~Smola}, {\em {Fastfood -- approximating kernel
  expansions in loglinear time}}, in Proceedings of the $30^{th}$ international
  conference on machine learning, JMLR: W\&CP volume 28, 2013.

\bibitem{lee}
{\sc J.~A. Lee and M.~Verleysen}, {\em Nonlinear dimensionality reduction},
  Springer, 2007.

\bibitem{mackay}
{\sc D.~J. MacKay}, {\em Introduction to gaussian processes}, NATO ASI Series F
  Computer and Systems Sciences, 168 (1998), pp.~133--166.

\bibitem{mauroy}
{\sc A.~Mauroy and I.~Mezi\'{c}}, {\em {On the use of Fourier averages to
  compute the global isochrons of (quasi) periodic dynamics}}, Chaos: An
  Interdisciplinary Journal of Nonlinear Science, 22 (2012), p.~033112.

\bibitem{mauroy2}
{\sc A.~Mauroy, I.~Mezi\'{c}, and J.~Moehlis}, {\em Isostables, isochrons, and
  koopman spectrum for the action-angle representation of stable fixed point
  dynamics}, Physica D: Nonlinear Phenomena, 261 (2013), pp.~19--30.

\bibitem{mezic2005}
{\sc I.~Mezi{\'c}}, {\em Spectral properties of dynamical systems, model
  reduction and decompositions}, Nonlinear Dynamics, 41 (2005), pp.~309--325.

\bibitem{mezic2013}
{\sc I.~Mezi{\'c}}, {\em Analysis of fluid flows via spectral properties of the
  koopman operator}, Annual Review of Fluid Mechanics, 45 (2013), pp.~357--378.

\bibitem{nadler2}
{\sc B.~Nadler, S.~Lafon, R.~R. Coifman, and I.~G. Kevrekidis}, {\em Diffusion
  maps, spectral clustering and eigenfunctions of fokker-planck operators}, in
  NIPS, 2005.

\bibitem{nadler}
{\sc B.~Nadler, S.~Lafon, R.~R. Coifman, and I.~G. Kevrekidis}, {\em Diffusion
  maps, spectral clustering and reaction coordinates of dynamical systems},
  Applied and Computational Harmonic Analysis, 21 (2006), pp.~113--127.

\bibitem{otto}
{\sc S.~E. Otto and C.~W. Rowley}, {\em Linearly-recurrent autoencoder networks
  for learning dynamics}.
\newblock arXiv:1712.01378, 2017.

\bibitem{Rahimi2007}
{\sc A.~Rahimi and B.~Recht}, {\em {Random features for large-scale kernel
  machines}}.
\newblock NIPS 20, 2007.

\bibitem{Rahimi2008}
{\sc A.~Rahimi and B.~Recht}, {\em {Weighted sums of random kitchen sinks:
  replacing minimization with randomization in learning}}.
\newblock NIPS 21, 2008.

\bibitem{rasmussen}
{\sc C.~E. Rasmussen}, {\em Gaussian processes in machine learning}, in
  Advanced lectures on machine learning, Springer, 2004, pp.~63--71.

\bibitem{Rohde}
{\sc C.~A. Rohde}, {\em {Generalized inverses of partitioned matrices}}, J.
  Soc. Indust. Appl. Math., 13 (1965), pp.~1033--35.

\bibitem{Rowley2009}
{\sc C.~W. Rowley, I.~Mezic, S.~Bagheri, P.~Schlatter, and D.~S. Henningson},
  {\em {Spectral analysis of nonlinear flows}}, J. Fluid Mech., 641 (2009),
  pp.~115--127.

\bibitem{SchmidJFM}
{\sc P.~J. Schmid}, {\em {Dynamic mode decomposition of numerical and
  experimental data}}, Journal of Fluid Mechanics, 656 (2010), pp.~5--28.

\bibitem{Schmidt2008}
{\sc P.~J. Schmid and J.~Sesterhenn}, {\em {Dynamic mode decomposition of
  numerical and experimental data}}, in 61$^{st}$ Annual Meeting of the APS
  Division of Fluid Dynamics, American Physical Society, 2008.

\bibitem{simonyan}
{\sc K.~Simonyan and A.~Zisserman}, {\em Very deep convolutional networks for
  large-scale image recognition}.
\newblock arXiv:1409.1556, 2014.

\bibitem{TuShangPIV_journal}
{\sc J.~H. Tu, C.~W. Rowley, J.~N. Kutz, and J.~K. Shang}, {\em Spectral
  analysis of fluid flows using sub-nyquist-rate piv data}, Experiments in
  Fluids, 55 (2014), p.~1805.

\bibitem{Tu2013_journal}
{\sc J.~H. Tu, C.~W. Rowley, D.~M. Luchtenburg, S.~L. Brunton, and J.~N. Kutz},
  {\em On dynamic mode decomposition: Theory and applications}, Journal of
  Computational Dynamics, 1 (2014), pp.~391--421.

\bibitem{Nystrom1}
{\sc C.~K.~I. Williams and M.~Seeger}, {\em {Using the Nystrom method to speed
  up kernel machines}}, Advances in Neural Information Processing Systems 13
  (NIPS 2000), MIT Press, 2001, pp.~682--88.

\bibitem{WilliamsEDMD_journal}
{\sc M.~O. Williams, I.~G. Kevrekidis, and C.~W. Rowley}, {\em A data--driven
  approximation of the koopman operator: Extending dynamic mode decomposition},
  Journal of Nonlinear Science, 25 (2015), pp.~1307--1346.

\bibitem{WilliamsKDMD_journal}
{\sc M.~O. Williams, C.~W. Rowley, and I.~G. Kevrekidis}, {\em A kernel-based
  method for data-driven koopman spectral analysis}, Journal of Computational
  Dynamics, 2 (2015), pp.~247--265.

\bibitem{wynn}
{\sc A.~Wynn, D.~S. Pearson, B.~Ganapathisubramani, and P.~J. Goulart}, {\em
  Optimal mode decomposition for unsteady flows}, Journal of Fluid Mechanics,
  733 (2013), pp.~473--503.

\bibitem{Nystrom2}
{\sc T.~Yang, Y.~F. Li, M.~Mahdavi, R.~Jin, and Z.~H. Zhou}, {\em {Nystrom
  method vs random Fourier features: a theoretical and empirical comparison}},
  Advances in Neural Information Processing Systems (NIPS), 2012, pp.~476--84.

\end{thebibliography}

\end{document}